\documentclass[sn-mathphys,Numbered]{sn-jnl}

\usepackage{graphicx}%
\usepackage{multirow}%
\usepackage{amsmath,amssymb,amsfonts}%
\usepackage{amsthm}%
\usepackage{mathrsfs}%
\usepackage[title]{appendix}%
\usepackage{xcolor}%
\usepackage{textcomp}%
\usepackage{manyfoot}%
\usepackage{booktabs}%
\usepackage{algorithm}%
\usepackage{algorithmicx}%
\usepackage{algpseudocode}%
\usepackage{listings}%
\usepackage{placeins}%

\usepackage{enumitem}
\usepackage{subcaption}
\usepackage{caption}
\usepackage{float}

\usepackage[most]{tcolorbox}

\begin{document}

\title[ASM for large SAAs of CCPs]{Active Set methods for solving large sample average approximations of chance constrained optimisation problems\footnote[2]{This paper sets forth the material presented at the \textit{Decision Making Under Uncertainty} workshop held at the University of Queensland, 11 July 2022.}}

\author*[1]{\fnm{Rick} \sur{Jeuken}}\email{f.jeuken@uq.net.au}

\author[1]{\fnm{Michael} \sur{Forbes}}\email{m.forbes@uq.edu.au}

\affil[1]{\orgdiv{School of Mathematics and Physics}, \orgname{University of Queensland}, \orgaddress{\street{St Lucia}, \city{Brisbane}, \postcode{4072}, \state{Queensland}, \country{Australia}}}

\abstract{This article describes a novel approach to chance-constrained programming based on the sample average approximation (SAA) method. Recent work focuses on heuristic approximations to the SAA problem and we introduce a novel approach which improves on some existing methods. Our Active Set method allows one to solve SAAs of chance-constrained programs with very large numbers of scenarios quickly. We demonstrate that increasing the number of scenarios is more important than improving accuracy with small numbers of scenarios.\\

We use an example of the portfolio selection problem to demonstrate the relative performance of previous and new methods. Extending the Active Set method to an integer-programming model further highlights its applicability and further improves over previous approaches.}


\maketitle

\section{Introduction}

Chance-constrained programming (CCP) has been in use since Charnes and Cooper \cite{CharnesCooper1959} published their original work. A chance constraint is one where random variables are constrained to a likelihood of being within a given tolerance. When stochasticity is known or estimated from data, chance constraints often provide a suitable avenue to characterise a problem. Situations with uncertain supply or demand, material properties, travel times etc. are all potential applications of CCP. Large volumes of data are being collected in many domains so it is increasingly easy to determine the distribution of uncertain data, or indeed to simply sample from historical data.\\ 

Adapted from \cite{Pagnoncelli2012}, a chance constrained optimization problem is generally of the form,\\
\\
CCP
\begin{equation}
min \; f(x) \label{CCP}
\end{equation}
\begin{equation*}
s.t. \:  Prob \left\lbrace G(x,r) \leq 0 \right\rbrace \geq 1 - \epsilon,
\end{equation*}
\begin{equation*}
x\in X,
\end{equation*}\\
\\
where $X \in \mathbb{R} ^n$ is a convex and closed set,  $r$ is a random vector with a probability distribution $P$ supported on a set $\Xi \subset \mathbb{R} ^d$, $\epsilon \in (0,1)$, $f:\mathbb{R} ^n  \to \mathbb{R}$ is a real valued function and $G: \mathbb{R} ^n \times \Xi \to \mathbb{R} ^m$.\\

For normally distributed random variables, writing the chance constraint as a quadratic constraint means that in some cases it is easily solvable in a modern solver, \cite{gurobi,Jeuken2021}. Many applications, however, involve random variables whose distribution is not normal or only known through sampling. Recent papers such as \cite{meng2020,stuhlmacher2020,wang2020,wu2020,zhang2021} describe some practical examples.\\ 

For random variables, not normally distributed, some methods transform the chance constraints to deterministic approximations using uniform or triangular distributions \cite{Zhang2019}. Some approaches formulate the problem as non-linear programs and solve those using heuristic methods such as genetic algorithms \cite{Chowdhury2019,Wang2018}. Several authors, \cite{bruninx2017, jiang2017,karimi2021,pourahmadi2019,stuhlmacher2020} apply the Sample Average Approximation (SAA) method to approximate chance-constrained programs with mixed-integer linear programs (MILPs). The SAA method is versatile and works for random variables of any distribution.\\

When formulating a chance constraint using the SAA method, a set of constraints replaces each chance constraint with one constraint for every sample of the random variables. The number of samples depends on the confidence and risk levels required. Formulating these problems is straightforward, however the increase in scenarios required for a given confidence or risk means these problems are NP-hard \cite{Luedtke2008}. Several authors focus on the minimum number of samples required to approximate the underlying problem to a sufficient level of accuracy. Campi and Garatti \cite{Campi2008}, Luedtke and Ahmed \cite{Luedtke2008} and Pagnonceill, Ahmed and Shapiro \cite{Pagnoncelli2009}, did much of the foundational work on quantifying the number of scenarios required.\\ 

Even with minimum numbers of sample scenarios, the size of the resultant problem can be overwhelming. Practical examples such as \cite{stuhlmacher2020} limit the application to small example models. Models that better approximate practical problems become unrealistic to solve as the scenario count increases. Recent approaches that further reduce the problem size look at methods to aggregate or cluster constraints \cite{chen2020}. Another approach is to discard constraints, usually by selecting those that are not binding to the solution, see \cite{Campi2011} and \cite{Kudela2020}. These discard methods are the current state of the art for general applicability.\\

Our work, outlined in this paper, shows that a simple modification of the methods in \cite{Pagnoncelli2012} and \cite{Kudela2020} can speed up solution times 7-fold. Introducing our Active Set method leads to further improvements in solution time and solution quality. We then extend the examples of \cite{Pagnoncelli2012} and \cite{Kudela2020} and the Active Set method to an integer-programming model. We thus demonstrate the versatility of the method. Computational experiments on an example of the Portfolio Selection Problem validate all methods presented.\\

\section{Previous Work}

The SAA method is a well understood and documented approximate approach to stochastic programming e.g \cite{shapiro2021}. The SAA method replaces a single probabilistic constraint with $N$ independent identically distributed scenarios. Each scenario is a distinct sample from the random part, i.e. the vector $r$ in CCP. We define the scenario program of CCP with $N$ scenarios as:\\
\\
SAA\textsubscript{N}
\begin{equation}
min \; f(x) \label{SAA_N}
\end{equation}
\begin{equation*}
s.t. \:  G(x,r^i) \leq 0, \;\;\; i=1,...,N,
\end{equation*}\\
\\
The SAA method is useful when there is no exact solution available, e.g. if the random variables do not have a “nice” Gaussian distribution. In the SAA method, a set of scenario constraints replaces each chance constraint. In early work, Campi and Garatti \cite{Campi2008} define a minimum number of scenarios that must be satisfied to ensure with a specified likelihood that the solution SAA\textsubscript{N} is feasible to the true problem. Pagnoncelli, Ahmed and Shapiro \cite{Pagnoncelli2009} go on to describe an example model of the portfolio problem and computational testing of the minimum scenario method.\\

Pagnoncelli, Ahmed and Shapiro \cite{Pagnoncelli2009} show it is a conservative approach to select $N$ scenarios and satisfy them all. Feasibility to the true problem is certain at a low risk, however the realised objective values are generally poor. Campi and Garatti \cite{Campi2011} extend the approach of replacing a chance constraint with a minimum number of scenario constraints. Instead of honouring all constraints, a select number can fail. They established that, if using $N$ scenario constraints and removing $k$, the solution satisfying the remaining constraints is feasible to the true problem CCP with confidence $(1-\beta)$ so long as $N$ and $k$ satisfy the condition:\\  

\begin{equation}
\binom{k+n-1}{k} \sum_{j=0}^{k+n+1} \binom{N}{j} \epsilon^j (1-\epsilon)^{N-j} \leq \beta \label{eqn3}
\end{equation}\\

Here $\beta \in (0,1)$ is any small confidence parameter value. Recall too that $n$ is the number of dimensions of $x$ in CCP and $\epsilon$ is the risk factor.\\ 

As $N$ goes to infinity we can choose $k/N$ arbitrarily close to $\epsilon$  with feasibility to the original problem. Pagnoncelli, Reich and Campi \cite{Pagnoncelli2012} show that when $n=20$,  $\beta = 10^{-9}$ and $\epsilon=0.05$ then over 1 million scenarios give $k/N$ greater than $0.045$. Perhaps more importantly, for a fixed value of $\beta$, as we increase $N$ and $k/N$ the expected objective value improves.\\ 

Using the relationship (\ref{eqn3}), the model SAA\textsubscript{N} is formulated as an integer program, using $N$ binary variables, $Z$:\\ 
\\
SAA\textsubscript{N,k}
\begin{equation}
min \; c^{T}x \label{SAA_Nk}
\end{equation}
\begin{equation*}
s.t. \:  G(x,r^i) \leq MZ_{i}, \;\;\; i=1,...,N,
\end{equation*}
\begin{equation*}
\sum_{i=1}^{N} Z_{i} \leq k
\end{equation*}
\begin{equation*}
Z_{i} \in \{0,1\}^{N}
\end{equation*}\\
\\
Here $M$ is a positive constant greater than the maximum possible value that $G(x,r)$ can attain.  The binary variables allow the failure of up to $k$ constraints. The solution of SAA\textsubscript{N,k} finds the best possible choice of the $k$ constraints to discard. As $k$ and $N$ satisfy (\ref{eqn3}) for the required $\epsilon$ we can be certain with a confidence $(1-\beta)$ that the solution is feasible to the true problem, CCP.\\

For the solution of SAA\textsubscript{N,k} one would naturally like the best objective value obtainable. The objective value improves as $N$ increases and this comes at the expense of tractability. Even if $G(x,r^{i})$ is linear, SAA\textsubscript{N,k} is still a Mixed Integer Programming (MIP) problem with $N$ binary variables. This quickly becomes unrealistic to solve. The idea is then, how to best approximate SAA\textsubscript{N,k} to obtain good results in a reasonable time. This approximation of an approximation is the approach of several previous authors and one that we extend in this paper.  Earlier work has mainly focussed on more accurate approximations of the solution of SAA\textsubscript{N,k} for a fixed $N$.  In this paper we focus on getting good approximate solutions for very large $N$.\\

Pagnoncelli, Reich, and Campi \cite{Pagnoncelli2012} provide detailed computational experiments of the constraint removal method applied to the portfolio selection problem. They trial two constraint removal schemes, one greedy and one randomised. We will henceforth refer to these as GR-P and RA-P for brevity. Both methods start with modelling SAA\textsubscript{N} as a linear program with $N$ scenario constraints.\\

The GR-P method iteratively removes $k$ constraints. At each iteration $i$, solving model SAA\textsubscript{N,k} determines the set of $n_{i}$ active constraints. Removing the active constraints one at a time and solving the corresponding model identifies the constraint whose removal best improves the objective value.  After permanently discarding the identified constraint, the entire process is repeated $k$ times. The GR-P algorithm requires solving $1+\sum_{i=1}^{k} n_{i}$  problems. In the RA-P algorithm, removing a random active constraint at each iteration reduces the number of problems to $1+k$.\\

Kudela and Popela \cite{Kudela2020} introduce a new and novel approach to the problem called the Pool and Discard algorithm. In their method, a reduced model has constraints added in rather than removed, which substantially reduces the solution times. Solving a similar number of models to the GR-P algorithm is required. Instead of each model having approximately $N$ constraints, each model starts with no constraints and adds them one at a time, until feasibility conditions are satisfied. Reducing the size of each model leads to the time improvement.\\ 

We note that the scenario removal methods derive from the theorem of Levin \cite{levin1969}. We can write the theorem of Levin as “The number of support scenarios for SAA\textsubscript{N} is at most $n$, the dimension of the decision vector $x$”. The trick is to find this minimal set of support scenarios. This paper contributes an efficient way of finding good approximate solutions to SAA\textsubscript{N,k} for very large $N$.\\

In this paper, we substantially improve the GR-P method and the Pool and Discard algorithm. We reduce the solution times by up to $85\%$. We also introduce a new algorithmic approach; an Active Set method that allows us to quickly obtain good approximate solutions for very large values of $N$. This further reduces the solution times for a comparable objective value. Computational experiments using an example of the portfolio selection problem demonstrate the advantages of the active set method. A further extension to solving an integer-programming variant of the portfolio problem highlights the versatility of our approach.\\

\section{Solution Algorithms}

\subsection{Active Set Method}

The explicit model for SAA\textsubscript{N} is very large with a constraint for each scenario. The GR-P and RA-P algorithms start with this large linear program and re-solve it repeatedly to determine which $k$ constraints are best to remove. The Pool and Discard Algorithm starts with no constraints and re-solves repeatedly to determine which constraints to add. The models are much smaller however the number of solves is similar to the number required for the GR-P and RA-P algorithms. The Active Set method reduces the overall work by starting with a limited model and only adding constraints. Both the size and the number of models will be small.\\

The Active Set method, like the Pool and Discard algorithm, starts with a relaxed version of SAA\textsubscript{N} that has no scenario constraints. We choose appropriate values of $N$, $k$ and $\epsilon$ for the instance of the problem and the required certainty level $\beta$ and these values satisfy relationship (\ref{eqn3}).\\

Solving the relaxed model, we test the solution against the scenarios and find a set of constraint violation values. If there are less than k constraint violations, then the incumbent solution satisfies SAA\textsubscript{N,k}. Otherwise we order this set of violation values from largest to smallest and select an element between the $k+1^{th}$ value and the end of the set of violations. The scenario of the selected value is added as a constraint to the set of constraints on  SAA\textsubscript{N}. The position of the selected value is a fixed proportion of the range from the $k+1^{th}$ value to the end of the set. We write the active set method explicitly as:\\ 

\begin{tcolorbox}[breakable,colback=lightgray,colframe=darkgray,width=\dimexpr\textwidth+12mm\relax,enlarge left by=-6mm] 

\textbf{Active Set Method}

\begin{enumerate}
\item Solve the relaxed version of SAA\textsubscript{N}:

\begin{equation}
min \; f(x) \label{5}
\end{equation}

and obtain a solution $\hat{x}$.

\item 	Test the solution of (\ref{5}) by computing the set of outcomes $O$, against each of the $N$ scenarios:

\begin{equation*}
O_{i} = G(x,r^{i}), \;\;\; i=1,…,N.
\end{equation*}
   
\item 	Order the set of outcomes where $O_{i}>0$ from largest to smallest, now $O^{R}$, while keeping track of the original scenario index for each outcome. 

\item 	If $|O^{R}| \leq k$, where $k$ is determined from (\ref{eqn3}) for the chosen $N$ and required $\epsilon$ and $\beta$, then stop.

\item  Otherwise: 

\begin{enumerate}

\item Determine $j$, where $j = \lfloor w(k+1) +(1-w) |O^{R} | \rfloor$ and $w \in [0,1]$.

\item Add the constraint $G(x,r^{p} )\leq 0$ to (\ref{5}), where $p$  is the original scenario index of $O_{j}^{R}$. Add the scenario index $p$ to a constraint set $C$. Re-solve and return to step 2. \\

\end{enumerate}
\end{enumerate}
\end{tcolorbox}

Of interest is how to decide which constraints to add. In steps 2-4 of the Active Set method we test the incumbent solution for feasibility to the true problem with certainty $(1-\beta)$. If the solution is feasible, we finish. If the solution is not feasible, we choose a violated scenario as a constraint.  If we consider a one-dimensional problem, then choosing the $k+1^{th}$ constraint eliminates all scenarios with worse outcomes. The $k+1^{th}$ cut does not eliminate the feasible scenarios in the single variable problem.\\

We extend the idea to $n$-dimensions, as a heuristic approach to selecting a minimal number of constraints. Now a violated constraint is selected that lies between the $k+1^{th}$ and the least violated scenarios. Varying the position of the selected scenarios (the value of the weight $w$) allows one to tune the method. Computational testing in Section 4 demonstrates the validity of the Active Set method.\\

The Active Set method described ends with a set of constraints, $C$. The final set of constraints in $C$ are not necessarily all binding. Iterating through and testing each constraint for removal or replacement may further improve the objective value. We can extend the Active Set method to check if removing or replacing some constraints further improves the solution. We describe the polishing step to the Active Set method: \\ 

\begin{tcolorbox}[breakable,colback=lightgray,colframe=darkgray,width=\dimexpr\textwidth+12mm\relax,enlarge left by=-6mm] 

\textbf{Active Set Method - Polish}\\

We start with the solution obtained in the Active Set method, $\hat{x}$ and the constraint set, $C$, i.e. satisfying:

\begin{equation}
min \; f(x) \label{6}
\end{equation}
\begin{equation*}
s.t. \:  G(x,r^i) \leq 0, \;\;\; i \in C.
\end{equation*}\\

For a fixed number of iterations:

\begin{enumerate}

\item 	For each constraint $s \in C$:
 
\begin{enumerate}

\item Remove constraint $s$ from (\ref{6}) and solve resulting model.

\item Test solution by computing the set of outcomes against each of the $N$ scenarios:

\begin{equation*}
O_{i} = G(x,r^i), \;\;\; i=1,...,N.
\end{equation*}
 
\item Order the set of outcomes from largest to smallest, now $O^{R}$.

\item Test the $k^{th}$ element of $O^{R}$ where $k$ is determined from (\ref{eqn3}) for the chosen $N$ and required $\epsilon$ and $\beta$. 

\item If $O_{k}^{R} > 0$ then: 

\begin{enumerate}

\item 	Replace $C_{s}$ with $G(x,r^{j}) \leq 0$, where $j$ is the original scenario index of $O_{k}^{R}$, as a constraint in (\ref{6}) and solve the resulting model. 

\item Test solution by computing the set of outcomes against each of the $N$ scenarios:

\begin{equation*}
O_{i}=G(x,r^{i}), \;\;\; i=1,...,N.   
\end{equation*}

\item Order the set of outcomes from largest to smallest, to create a new $O^{R}$.

\end{enumerate}

\item If $O_{k}^{R} \leq 0$ and the new solution is better than the incumbent solution, $\hat{x}$ is the new incumbent.\\ 

\end{enumerate}
\end{enumerate}
\end{tcolorbox}

\subsection{Dual Improvements}

Some of the historical constraint removal schemes suffer from the curse of dimensionality. The original GR-P algorithm operates on a model with $N$ constraints and requires approximately $1+ \sum_{i=1}^{k} n_{i}$ iterations of solving. The original work in \cite{Pagnoncelli2012} shows that as $N$ increases then $k/N$ approaches $\epsilon$. Thus, the proportion of allowed failures increases and we obtain a better objective value. As $N$ increases however, the solution times of the GR-P algorithm increase. At some point, it is unrealistic to continue. If solution times are of the order of hours or days the algorithms usefulness is limited.\\

The Pool and Discard method improves greatly over the GR-P and RA-P algorithms even though the number of model solves is similar. The number of constraints in each model is approximately $n$, much less than $N$. Reducing the model size substantially decreases the solution times.\\

When the underlying model is a linear program we can use dual variables to greatly speed up these historical methods.  In both models, there is an internal iterative loop removing binding constraints and re-solving. The removal that results in the best improvement of the objective value is selected. As an alternative to removing a binding constraint and resolving the model, we can use the dual variable associated with each binding constraint as an approximation of the improvement in the objective from removing that constraint.\\

In the GR-P and Pool and Discard algorithms we replace the internal loop with a set of ordered dual variables of the model constraints. We then remove the constraint with the dual value that best improves the objective value. The dual improvement only requires the solution of $1+k$ models. This requires the same number of model solves as RA-P and the random versions of the Pool and Discard algorithm. Instead of a random removal though (from the binding constraints), we select the constraint that gives the best instantaneous improvement to the objective value.\\

An alternative polishing phase in the Active Set method also uses the information we have in the dual values of each of the constraints in set $C$. Instead of removing each constraint and re-solving the model to obtain the best improvement to the objective value, we order the set of constraints using their dual values.\\

\begin{tcolorbox}[breakable,colback=lightgray,colframe=darkgray,width=\dimexpr\textwidth+12mm\relax,enlarge left by=-6mm] 

\textbf{Active Set Method - Polish using dual variables}\\

We again start with the solution obtained in the Active Set method, $\hat{x}$ and the constraint set, $C$, i.e. satisfying (\ref{6}):\\

For a fixed number of iterations:

\begin{enumerate}

\item Obtain the dual value for each constraint $s \in C$, i.e. $\pi_{s}$.
 
\item Find the constraint whose dual value indicates that removing this constraint has the greatest instantaneous improvement in the objective value.

\item Remove the constraint that has the dual value found in step 2 and solve the resulting model.  

\item Test solution by computing the set of outcomes against each of the $N$ scenarios:

\begin{equation*}
O_{i} = G(x,r^i), \;\;\; i=1,...,N.
\end{equation*}
 
\item Order the set of outcomes from largest to smallest, now $O^{R}$.

\item Test the $k^{th}$ element of $O^{R}$ where $k$ is determined from (\ref{eqn3}) for the chosen $N$ and required $\epsilon$ and $\beta$. 

\item If $O_{k}^{R} > 0$ then: 

\begin{enumerate}

\item 	Replace $C_{s}$ with $G(x,r^{j}) \leq 0$, where $j$ is the original scenario index of $O_{k}^{R}$, as a constraint in (\ref{6}) and solve the resulting model. 

\item Test solution by computing the set of outcomes against each of the $N$ scenarios:

\begin{equation*}
O_{i}=G(x,r^{i}), \;\;\; i=1,...,N.   
\end{equation*}

\item Order the set of outcomes from largest to smallest, to create a new $O^{R}$.

\end{enumerate}

\item If $O_{k}^{R} \leq 0$ and the new solution is better than the incumbent solution, $\hat{x}$ is the new incumbent.\\ 

\end{enumerate}

\end{tcolorbox}

The computational testing in Section 4 demonstrates the dual improvements we have made to historical methods. Comparing with our Active Set method, with and without the two polishing steps, gives useful insights to all methods.   

\subsection{Integer Models}

All the methods and models described so far, in historical work and our own Active Set method apply when the underlying optimisation problem is a linear program. However, if the underlying problem is an integer program then the solution times using methods of adding and removing constraints escalate rapidly. When we apply the methods discussed to an integer-programming model the advantage lies with an algorithm that requires fewer solves and that removes fewer constraints.\\

Similar sized integer programs generally take longer to solve than linear programs so a smaller model is better. Computational testing in Section 4 demonstrates applying the GR-P, RA-P and Pool and Discard algorithms, as well as our Active Set method, to an integer-programming model.\\

\section{Computational Models}

In this section, we compare the methods discussed in sections two and three. The computational example we use is a version of the Portfolio Selection Problem. 

\subsection{Portfolio Selection Problem}

Portfolio selection problems originated with Markowitz \cite{Markowitz1952} and have a long and storied history of application and approximation. There is an excellent review of previous material in \cite{Xidonas2020}. \cite{Pagnoncelli2009} and \cite{Pagnoncelli2012} also consider the portfolio problem. A formulation of the chance constrained portfolio selection, or asset allocation, problem is:\\
\\
AAP
\begin{equation}
\max\limits_{x \in X} \; \mathbb{E} [r^{T}x] \label{AAP}
\end{equation}
\begin{equation*}
s.t. \:  Prob\{r^{T}x \geq \alpha \} \geq 1 - \epsilon,
\end{equation*}
\begin{equation*}
\sum_{x \in X} x=1
\end{equation*}\\
\\

The vector $x$ represents the fraction of total wealth invested in each of $n$ available assets. The vector of random returns of each of the assets is $r$. The objective is to maximise the expected return subject to that return being greater than or equal to a desired $\alpha$, and for this to occur with certainty $(1-\epsilon)$. In practice the expected value of $r$ is slightly above 1 as each asset class increases in value over the investment period. The maximum acceptable loss is $\alpha$, which is often slightly below one. That is, we require the probability the portfolio reduces in value to below $\alpha$ to be no more than $\epsilon$. This is a simple formulation of the portfolio selection problem and we use it as an example to compare the methods and improvements described in this paper.\\ 

\subsection{Computational Trials}

We use AAP as an example to outline the introduced solution algorithms. We can approximate the true problem AAP by sampling $N$ scenarios from the random vector of returns $r$ and adding a constraint for each scenario. Allowing $k$ constraints to fail satisfies (\ref{eqn3}) for a given $\epsilon$ and $\beta$. Writing AAP as an integer-programming model, with $M$ as a positive constant greater than the largest possible failure, and $N$ binary variables $(Z$):\\  
\\
SAA\textsubscript{N,k}\textsuperscript{P}
\begin{equation}
\max \; \mathbb{E} [r^{T}x] \label{SAA_Nk^P}
\end{equation}
\begin{equation*}
s.t. \sum r_{s}^{T} x \geq \alpha - M Z_{s} \;\;\; s=1,...,N,
\end{equation*}
\begin{equation*}
\sum_{s=1}^{N} Z_{s} \leq k,
\end{equation*}
\begin{equation*}
\sum_{i=1}^{n} X_{s} = 1,
\end{equation*}
\begin{equation*}
Z_{s} \in \{0,1\}^{N}
\end{equation*}\\
\\

Solving SAA\textsubscript{N,k}\textsuperscript{P} exactly allows us to determine the best set of $k$ constraints we can allow to fail. Thus the solution of SAA\textsubscript{N,k}\textsuperscript{P} best approximates AAP. This model quickly becomes intractable however, as it requires $N$ binary variables.\\

The numerical experiments use twenty assets plus one risk free option, $(n=21)$. We obtained data from Yahoo\textsuperscript{TM} Finance for the twenty assets in the ASX-20 in early 2021. For all assets, we use monthly data from the earliest available year to early 2021. We take the price of each month over the price of one year prior. These quotients are the returns and we assume they follow a multivariate normal distribution. We estimate the parameters using unbiased estimators for means and covariance.\\

The risk-free asset is included to ensure that any model is feasible, so long as $\alpha \leq 1$. Assuming the distribution of the returns is multivariate normal is not strictly true. However, a multivariate normal model allows us to solve the true problem AAP directly as a second order conical program, SOCP. The exact solution is a benchmark to compare the various methods.\\
\\
SOCP
\begin{equation}
\max \sum_{i=1}^{n} X_{i} \bar{u_{i}} \label{SOCP}
\end{equation}
\begin{equation*}
\sum_{i=1}^{n} X_{i} = 1
\end{equation*}
\begin{equation*}F^{-1} (1-\epsilon)^{2} X^{T} Cov X \leq \left( \sum_{i=1}^{n} X_{i} \bar{u_{i}} - \alpha \sum_{i=1}^{n} X_{i} \right)^{2}
\end{equation*}\\
\\

Here $\bar{u_{i}}$ is the mean expected return of asset $i$, $Cov$ is the covariance matrix of the $n$ assets, $X$ is the decision vector of allocating assets and $F^{-1}$ is the inverse normal distribution.\\

We solve AAP using each of these methods:\\

\begin{enumerate}

\item SOCP, solved to obtain the exact solution. We compare the remaining methods with this baseline.

\item SAA\textsubscript{N,k}\textsuperscript{P}, the integer-programming version. This solution gives the best possible set of scenarios to remove.

\item GR-P, the greedy algorithm as described in \cite{Pagnoncelli2012}.

\item FGR-P, the greedy algorithm GR-P with the speed up that uses dual variables in the inner loop.

\item PND, the Pool and Discard algorithm from \cite{Kudela2020}.

\item FPND, the Pool and Discard algorithm with the speed up using dual variables in the inner loop.

\item ASM-1, our active set method described in section 3.1 without the polishing step. We set $w=0.5$ in step 5(a).

\item ASM-2, our active set method with the first polishing algorithm, using the solution of ASM-1 as the initial solution.

\item ASM-3, our active set method with the second polishing algorithm, using the dual variables to decide constraint removal and replacement. The solution of ASM-1 is the initial solution.\\

\end{enumerate}

For each of these nine methodologies, we test over 30 runs using the same set of scenarios for each trial. AAP is set up with $\alpha=0.05$ and $\epsilon=0.05$. For the scenarios we range values of  $N$ from 1000 to 1,000,000 and choose $k$ to obtain the maximum value of $\beta$ that is less than 0.000005.  These values are given in Table 1. The allowable time for each run is limited to one hour.\\

In each trial where more than one model is solved, all models after the first are warm-started. Each iteration involves adding or removing one of many constraints. Using the solution of the $i^{th}$ iteration as a starting basis for the $i^{th}+1$ iteration means the average solution time is reduced relative to starting each model from scratch. \cite{Kudela2020} provides evidence of the reduction.\\

\begin{table}[h] 

\caption{\label{tab:Scenarios}Values of $N$ and $k$ satisfying (\ref{eqn3}) that are less than or equal to the $\beta$ likelihood. Note that $\epsilon=0.05$ in (\ref{eqn3}).}

\begin{tabular}{|c|c|c|c|} 
\hline
\textbf{$N$} & \textbf{$k$} & \textbf{$\beta$} & \textbf{$k/N$}  \\ 
\hline
1000       & 1           & 1.53E-05\textsuperscript{(1)}  & 0.001          \\
2500       & 24          & 3.73E-06                       & 0.010          \\
5000       & 85          & 3.46E-06                       & 0.017          \\
10000      & 238         & 3.31E-06                       & 0.024          \\
20000      & 593         & 4.40E-06                       & 0.030          \\
50000      & 1786        & 3.75E-06                       & 0.036          \\
100000     & 3923        & 4.72E-06                       & 0.039          \\
500000     & 22278       & 4.96E-06                       & 0.045          \\
1000000    & 45978       & 4.74E-06                       & 0.046          \\
\hline
\end{tabular}

\footnotesize{(1) Note that $\beta$ for $N=1000$ is much larger than for $N \geq 2500$ because of the reduced flexibility in choosing $k$. If $k=0$, then $\beta \approx 2.88E-07$.}\\

\end{table}
In the SOCP models we use $k/N$ as an approximation of $\epsilon$. We do this to give the exact solution at that failure rate, to compare to the other methods. Using $\epsilon = 0.05$ in SOCP gives the exact solution to the true problem. We however, want the exact solution for the specific failure rate we are testing.\\

To test the quality of the solutions, we measure the probability of violation by evaluating each solution against a test set of 100,000 scenarios. The probability of violation is the proportion of these scenarios that fail, with a returned objective value less than $\alpha$. We also determine a binomial confidence interval to determine the $\beta$ likelihood upper limit of the interval, \cite{wallis2013binomial}. By definition this should always be $\leq \beta$.\\

We implement and solve all models on the same machine, an Intel i7-6700 CPU @ 3.4GHz with 8 threads and 16GB of RAM. The scenarios were evaluated using Gurobi 9.1.2 via Python scripting. When applicable, all MIP problems are solved to default optimality (tolerance set to 1e-4). \\
  
\subsection{Results}

We report the mean objective values, solution times, count of LP’s solved, probability of violation and the upper limit of the binomial confidence interval for each method and number of scenarios, Tables 2(a) – 2(e).\\

\vspace{-0.2in}

\begin{table}[h] 

\caption{Objective values, solution times, model solves probability of violation and confidence of violation for the methods, SOCP, SAA\textsubscript{N,k}\textsuperscript{P}, GR-P, FGR-P, PND, FPND, ASM-1, ASM-2 and ASM-3.}

\begin{subtable}[h]{1.0\textwidth} 
\begin{tabular}{|l|ccccccccc|} 
\hline
$N$       & SOCP   & SAA\textsubscript{N,k}\textsuperscript{P} & GR-P   & FGR-P  & PND    & FPND   & ASM-1  & ASM-2  & ASM-3   \\ 
\hline
1000    & 1.0302 & 1.0442                    & 1.0395 & 1.0392 & 1.0476 & 1.0415 & 1.0388 & 1.0392 & 1.0389  \\
2500    & 1.0494 & 1.0528                    & 1.0522 & 1.0518 & 1.0526 & 1.0526 & 1.0522 & 1.0530 & 1.0523  \\
5000    & 1.0608 & -                         & 1.0617 & 1.0612 & 1.0621 & 1.0613 & 1.0594 & 1.0612 & 1.0597  \\
10000   & 1.0712 & -                         & 1.0719 & 1.0716 & 1.0720 & 1.0722 & 1.0672 & 1.0696 & 1.0680  \\
20000   & 1.0805 & -                         & 1.0807 & 1.0803 & 1.0808 & 1.0805 & 1.0745 & 1.0773 & 1.0749  \\
50000   & 1.0910 & -                         & -      & 1.0903 & -      & 1.0903 & 1.0814 & 1.0848 & 1.0825  \\
100000  & 1.0976 & -                         & -      & -      & -      & -      & 1.0877 & 1.0894 & 1.0878  \\
500000  & 1.1084 & -                         & -      & -      & -      & -      & 1.0935 & 1.0971 & 1.0942  \\
1000000 & 1.1115 & -                         & -      & -      & -      & -      & 1.0971 & 1.1003 & 1.0977  \\
\hline
\end{tabular}

\caption{Objective values for the methods, SOCP, SAA\textsubscript{N,k}\textsuperscript{P}, GR-P, FGR-P, PND, FPND, ASM-1, ASM-2 and ASM-3. The exact solution to the true problem, solved using SOCP, returns an objective value of 1.1209.In each instance of the SOCP models we solve for $k/N$ in place of $\epsilon$. We do this to give the exact solution at that failure rate.}

\end{subtable}

\end{table}

\begin{table}\ContinuedFloat

\begin{subtable}[h]{1.0\textwidth} 

\begin{tabular}{|l|ccccccccc|} 
\hline
$N$ & SOCP & SAA\textsubscript{N,k}\textsuperscript{P} & GR-P   & FGR-P  & PND    & FPND   & ASM-1 & ASM-2  & ASM-3  \\ 
\hline
1000           & 13.4 & 0.6                       & 0.1    & 0.1    & 0.2    & 0.1    & 0.1   & 1.0    & 0.4    \\
2500           & 21.6 & 778.0                      & 2.8    & 0.6    & 6.4    & 1.2    & 0.3   & 3.2    & 0.9    \\
5000           & 27.3 & >3600                      & 20.2   & 3.7    & 42.4   & 6.7    & 0.6   & 7.5    & 1.7    \\
10000          & 31.7 & -                         & 129.8  & 20.8   & 235.2  & 35.5   & 1.5   & 18.1   & 3.8    \\
20000          & 31.4 & -                         & 1010.7 & 143.7  & 1163.8 & 169.0   & 3.4   & 40.4   & 8.2    \\
50000          & 31.6 & -                         & >3600   & 1149.1 & >3600   & 1258.2 & 10.6  & 119    & 23.3   \\
100000         & 32.1 & -                         & -      & >3600   & -      & >3600   & 23.9  & 241.2  & 51.2   \\
500000         & 32.0 & -                         & -      & -      & -      & -      & 118.6 & 1510.3 & 249.8  \\
1000000        & 32.3 & -                         & -      & -      & -      & -      & 246.8 & 3033.3 & 522.6  \\
\hline
\end{tabular}

\caption{\emph{Overall solution times for the methods, SOCP, SAA\textsubscript{N,k}\textsuperscript{P}, GR-P, FGR-P, PND, FPND, ASM-1, ASM-2 and ASM-3.}}

\end{subtable}

\begin{subtable}[h]{1.0\textwidth} 

\begin{tabular}{|l|ccccccccc|} 
\hline
$N$ & SOCP & SAA\textsubscript{N,k}\textsuperscript{P} & GR-P   & FGR-P & PND    & FPND   & ASM-1 & ASM-2 & ASM-3  \\ 
\hline
1000           & 1    & 1                         & 8.8    & 2     & 26.7   & 12.1   & 37.1  & 194.9 & 77.1   \\
2500           & 1    & 1                         & 211.1  & 25    & 355.8  & 65.5   & 45.6  & 258.6 & 85.6   \\
5000           & 1    & -                         & 786.9  & 86    & 1190.2 & 188.5  & 53.9  & 297.5 & 93.6   \\
10000          & 1    & -                         & 2249.1 & 239   & 3314.3 & 497.0  & 62.2  & 335.5 & 101.4  \\
20000          & 1    & -                         & 5733.1 & 594   & 8306.2 & 1208.0 & 71.3  & 376.7 & 110.6  \\
50000          & 1    & -                         & -      & 1787  & -      & 3594.4 & 83.5  & 395.4 & 121.6  \\
100000         & 1    & -                         & -      & -     & -      & -      & 96.7  & 424.7 & 136.3  \\
500000         & 1    & -                         & -      & -     & -      & -      & 120.0 & 503.1 & 156.5  \\
1000000        & 1    & -                         & -      & -     & -      & -      & 127.2 & 495.6 & 163.9  \\
\hline
\end{tabular}

\caption{\emph{Number of times a model is solved  the methods, SOCP, SAA\textsubscript{N,k}\textsuperscript{P}, GR-P, FGR-P, PND, FPND, ASM-1, ASM-2 and ASM-3.}}

\end{subtable}

\begin{subtable}[h]{1.0\textwidth}

\begin{tabular}{|l|cccccccccc|} 
\hline
$N$       & $k/N$ & SOCP  & SAA\textsubscript{N,k}\textsuperscript{P} & GR-P  & FGR-P & PND   & FPND  & ASM-1 & ASM-2 & ASM-3  \\ 
\hline
1000    & 0.001               & 0.001 & 0.010                      & 0.007 & 0.007 & 0.029 & 0.013 & 0.007 & 0.007 & 0.006  \\
2500    & 0.010                & 0.010  & 0.015                     & 0.014 & 0.013 & 0.015 & 0.015 & 0.014 & 0.015 & 0.014  \\
5000    & 0.017               & 0.018 & -                         & 0.020  & 0.019 & 0.021 & 0.020  & 0.020  & 0.020  & 0.02   \\
10000   & 0.024               & 0.025 & -                         & 0.026 & 0.026 & 0.026 & 0.027 & 0.026 & 0.026 & 0.026  \\
20000   & 0.030                & 0.030  & -                         & 0.031 & 0.031 & 0.031 & 0.031 & 0.031 & 0.031 & 0.031  \\
50000   & 0.036               & 0.036 & -                         & -     & 0.036 & -     & 0.036 & 0.036 & 0.036 & 0.036  \\
100000  & 0.039               & 0.040  & -                         & -     & -     & -     & -     & 0.040  & 0.039 & 0.04   \\
500000  & 0.045               & 0.045 & -                         & -     & -     & -     & -     & 0.045 & 0.044 & 0.045  \\
1000000 & 0.046               & 0.047 & -                         & -     & -     & -     & -     & 0.046 & 0.045 & 0.046  \\
\hline
\end{tabular}

\caption{\emph{Probability of Violation for the solution of each of the methods, SOCP, SAA\textsubscript{N,k}\textsuperscript{P}, GR-P, FGR-P, PND, FPND, ASM-1, ASM-2 and ASM-3.}}

\end{subtable}

\begin{subtable}[H]{1.0\textwidth}

\begin{tabular}{|l|cccccccccc|} 
\hline
$N$ & $k/N$  & SOCP  & SAAN,k\textsuperscript{P} & GR-P  & FGR-P & PND   & FPND  & ASM-1 & ASM-2 & ASM-3  \\ 
\hline
1000         & 0.001 & 0.001 & 0.012                     & 0.008 & 0.008 & 0.032 & 0.015 & 0.008 & 0.008 & 0.008  \\
2500         & 0.010  & 0.011 & 0.017                     & 0.015 & 0.015 & 0.017 & 0.016 & 0.016 & 0.016 & 0.016  \\
5000         & 0.017 & 0.020  & -                         & 0.022 & 0.021 & 0.023 & 0.022 & 0.022 & 0.022 & 0.022  \\
10000        & 0.024 & 0.027 & -                         & 0.028 & 0.028 & 0.028 & 0.029 & 0.029 & 0.028 & 0.029  \\
20000        & 0.030  & 0.033 & -                         & 0.033 & 0.033 & 0.033 & 0.033 & 0.034 & 0.034 & 0.034  \\
50000        & 0.036 & 0.039 & -                         & -     & 0.039 & -     & 0.039 & 0.039 & 0.039 & 0.039  \\
100000       & 0.039 & 0.043 & -                         & -     & -     & -     & -     & 0.043 & 0.042 & 0.043  \\
500000       & 0.045 & 0.048 & -                         & -     & -     & -     & -     & 0.048 & 0.047 & 0.048  \\
1000000      & 0.046 & 0.050  & -                         & -     & -     & -     & -     & 0.049 & 0.048 & 0.049  \\
\hline
\end{tabular}
\caption{\emph{Confidence of violation (upper limit) for the solution of each of the methods, SOCP, SAA\textsubscript{N,k}\textsuperscript{P}, GR-P, FGR-P, PND, FPND, ASM-1, ASM-2 and ASM-3.}}
\end{subtable}

\end{table}

\FloatBarrier

Some immediate observations are clear when looking at the results. The SOCP model provides the benchmark “true” or “best” result. The SAA\textsubscript{N,k}\textsuperscript{P} model is likely a good one but completely intractable when considering the solution times. The remaining heuristic methods have a common thread. The longer it takes, the better the solution, this appears inescapable.\\

The average number of solutions required, Table 2(c), clearly defines the reason for overall run times – algorithms requiring many solves take longer. One can also see, Table 2(d), that testing the scenarios gives a performance that relates closely to $k/N$ as expected from the definition of the method. Then one sees in Table 2(e) the upper limit of the binomial confidence interval for the $\beta$ likelihood. Even for very large $N$ we maintain a solution to the true problem with greater than $(1-\beta)$ likelihood. Thus we validate all of the methods.\\   

What is clear Figure 1, is that our introduced active set method delivers the best solution in the fastest time. While the other approximate methods usually give a better answer in the end, this often takes much longer. In any time/value trade off, using our active set methods delivers the better results. When considering problems requiring integer programming, the time/value trade-off is more evident.\\

\begin{figure}[h] 
\begin{center}
\includegraphics[scale=0.7]{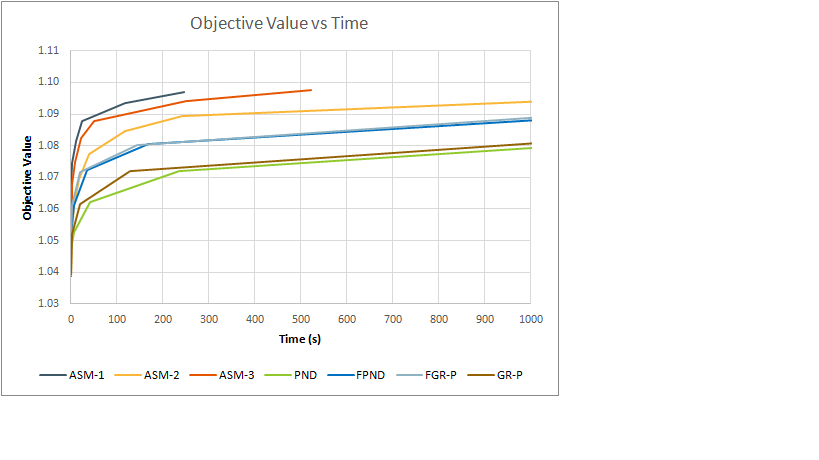}
\end{center}
\vspace{-0.2 in}
\caption{Plot of Objective Value vs total solution time for the heuristic methods GR-P, FGR-P, PND, FPND, ASM-1, ASM-2 and ASM-3. }
\end{figure}

Selecting the weight value in step 5(a) of the Active Set method has an impact on the outcome of running the algorithms. Figure 2 shows some different selections of $w$ for the case of $N=100,000$. One sees some trade-off where a better objective value is obtained a little more slowly by selecting the scenario constraint at the midpoint between the $k+1^{th}$ violated and least violated scenarios. Selecting a position close to the least violated, $w=0.01$, takes a long time as many scenario constraints are added. Selecting the $k+1^{th}$ most violated gives a lower objective value, on average, as too many feasible solutions are cut off.\\

\begin{figure}[h]
\begin{center}
\includegraphics[scale=0.7]{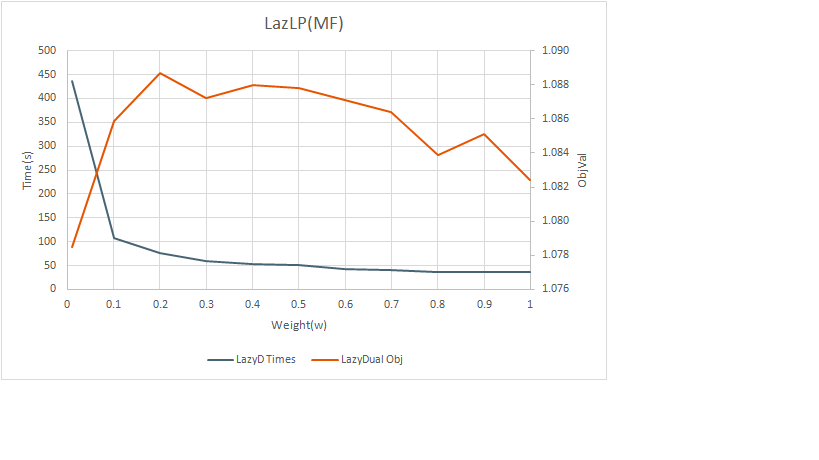}
\end{center}
\vspace{-0.5 in}
\caption{Objective values, solution times and number of constraints added for different values of w in the method ASM-1, with $N = 100,000$, $k=3923$. }
\end{figure}

\subsection{Remarks on Complexity}

We observe in Table 2(b) that the time to solve the example problem is always considerably less using the Active Set method. We introduce the idea that the solution times are proportional to the number of linear programming (LP) models each algorithm solves. These simple linear programs solve quickly in solvers and there is some relatively fixed amount of computational effort that sits around each solution of a LP.\\

There are several common features in the methods. When k scenarios are being removed from a constraint set, there are $k$ iterations required. Where binding constraints are being tested for removal, the best case is one constraint to test and possibly remove and the worst case is $n$ constraints, \cite{levin1969}. In Table 3 we see the best and worst case number of iterations requiring the solution of a linear program.\\     

\begin{table}[h] 

\caption{Computational complexity, best and worst case for the methods, GR-P, FGR-P, PND, FPND, ASM-1, ASM-2 and ASM-3. For ASM-2 and ASM-3, $a$ is the fixed number of iterations selected.}

\begin{tabular}{|l|cc|} 
\hline
Method & Best Case & Worst Case    \\ 
\hline
GR-P   & $k+1$      & $kn+k+1$       \\
FGR-P  & $k+1$      & $k+1$          \\
PND    & $k$        & $k(n+n^{2})$  \\
FPND   & $k$        & $kn$           \\
ASM-1  & $2$        & $N-k+1$        \\
ASM-2  & $2+a$      & $N-k+1+2an$    \\
ASM-3  & $2+a$      & $N-k+1+2a$     \\
\hline
\end{tabular}
\end{table}

In Table 3 it is evident the best case for the historical algorithms is approximately $k$ solutions required whereas the best case for ASM-1 is only two iterations. ASM-2 and ASM-3 add some small selected value for $a$. The best case for ASM-1 requires solving a model without constraints, adding a constraint, testing and accepting the solution. We plot the solve count vs the scenario count for the different algorithms applied to the example problem in Figure 3. For the methods GR-P, FGR-P, PND and FPND there is an observed linear increase in complexity as the scenario count $N$ increases. We know it will always be greater than $k$ and this increases as $N$ increases, up to a limit of $\epsilon N$. The observed empirical performance of the active set methods is $log(N)$.  

\begin{figure}[H] 
\centering
\begin{subfigure}{0.45\textwidth}
\includegraphics[scale=0.4]{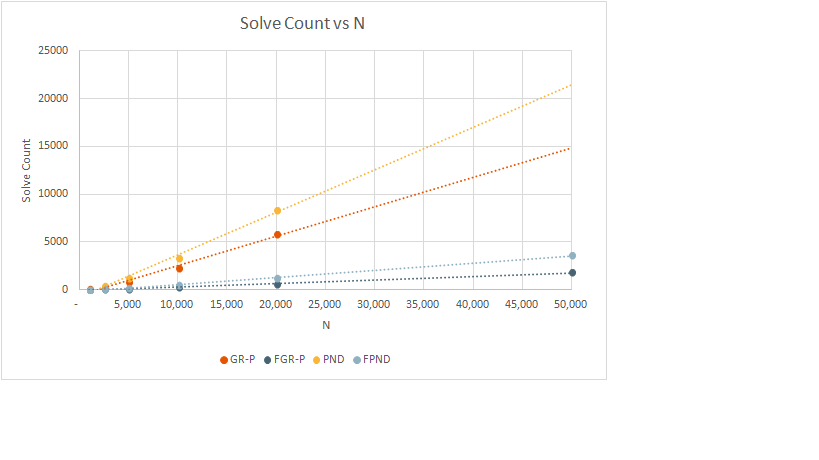}
\caption{\emph{Plots for GR-P, FGR-P, PND, FPND. Note the linear scale of the x-axis}}
\end{subfigure}
\hfill
\begin{subfigure}{0.45\textwidth}
\includegraphics[scale=0.4]{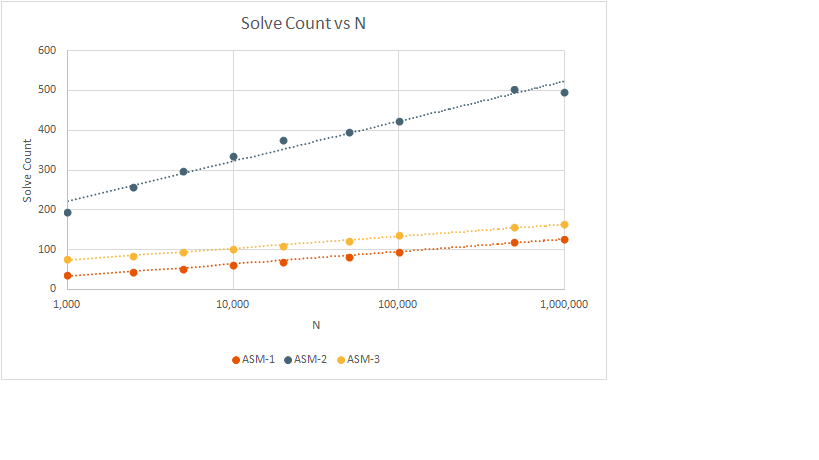}
\caption{\emph{Plots for ASM-1, ASM-2 and ASM-3. Note the log scale of the x-axis}}
\end{subfigure}
\caption{Plots of the Solve Count vs $N$.}
\vspace{-0.2 in}
\end{figure}

All of the described methods are trying to find the ‘best’ set of no more than n binding constraints. With this ‘best’ set, any additional constraints may be discarded without improving the solution and the removal of one will result in an infeasible solution, \cite{levin1969}. Methods which rely on selecting and adding or removing $k$ constraints have a complexity that is necessarily greater than $k$. A method such as our Delayed Constraint method looks for a set of approximately $n$ constraints and this is efficient to solve problems such as the portfolio selection problem examined in this paper.\\ 

Our Delayed Constraint method has an observed number of model solutions of order $log(N)$ and this remains readily solvable for very large values of $N$. We also see, in Table 2a, that solving problems with very large values of $N$ gives solutions with better objective values than more precise methods that can only handle relatively small values of $N$.\\

\subsection{An Integer-programming Model}

All of the described methods are trying to find the ‘best’ set of no more than n binding constraints. With this ‘best’ set, any additional constraints may be discarded without improving the solution and the removal of one will result in an infeasible solution \cite{levin1969}. Methods which rely on selecting and adding or removing $k$ constraints have a complexity that is necessarily greater than $k$. A method such as our Active Set method looks for a set of approximately $n$ constraints and this is efficient to solve problems such as the portfolio selection problem examined in this paper.\\ 

Our Active Set method has an observed number of model solutions of order $log(N)$ and this remains readily solvable for very large values of $N$. We also see, in Table 2a, that solving problems with very large values of $N$ gives solutions with better objective values than more precise methods that can only handle relatively small values of $N$.\\

\subsection{An Integer-programming Model}

Six of the nine methods discussed are directly applicable if the problem presents as an integer-programming model. The other three methods (FGR-P, FPND and ASM-3) use dual variables and so are not applicable.  The SOCP model is still available if the random vector is multivariate normally distributed. The SAA\textsubscript{N,k}\textsuperscript{P} method is essentially unchanged, with some added constraints. Applying the heuristic methods, we see some structural changes to the solution times.\\

We update the model of the Portfolio Problem, AAP, by adding a simple constraint:\\

\begin{equation}
X_{i} \; \in \{ 0,[l,u] \} \;\;\; \forall i=1,…,n-1. \label{AAP-Int}
\end{equation}\\
     
Constraint (\ref{AAP-Int}) requires, if investing in an asset, then that investment must be greater than $l$ and less than $u$ where $ 0<l<u<1$. The added constraint reflects that a decision maker may want to remove many small investments and limit the exposure to any single asset. The $n^{th}$ asset is the cash option, and this stays unconstrained to ensure the problem is always feasible. We implement constraint (\ref{AAP-Int}) with $n-1$ binary variables, thus modelling the problem as an integer program.\\

To demonstrate the relative performance of the methods already described, we solve AAP with constraint (\ref{AAP-Int}), where practical, by:

\begin{enumerate}

\item Adding constraint (\ref{AAP-Int}) to SOCP to obtain the exact solution. 
\item Adding constraint (\ref{AAP-Int}) to SAA\textsubscript{N,k}\textsuperscript{P} to find the best set of constraints to remove
\item GR-P
\item PND 
\item ASM-1, using a weight, $w=0.5$. 

\end{enumerate}

For each method, we again evaluate 30 trials using the same sets of scenarios as the LP trials (Table 3). Several of the heuristic methods require solving a model many times. The greedy methods GR-P and ASM-2 require re-solving with each removal of a binding constraint. When the model is an LP, each resolve uses warm starts and is relatively quick. The algorithm then runs through in reasonable time (Table 2(b)). When the model is an IP every solve takes longer and the overall time becomes unrealistic to implement. We do not present the trials of ASM-2 for the IP due to solution time violations. There are incomplete result sets for GR-P and PND for the same reason.\\

The results using the described models and methods SOCP, SAA\textsubscript{N,k}\textsuperscript{P}, GR-P, PND and ASM-1 to obtain solutions to AAP with constraint (\ref{AAP-Int}) are in Tables 4a – 4f. Run times are capped - we do not evaluate models with solution times greater than one hour. All values given are the mean of the 30 trials.\\

\begin{table}[h] 

\caption{Objective values, solution times, model solves probability of violation and confidence of violation for the methods, SOCP, SAA\textsubscript{N,k}\textsuperscript{P}, GR-P, PND and ASM-1.}

\begin{subtable}[h]{1.0\textwidth}

\begin{tabular}{|l|ccccc|} 
\hline
$N$       & SAA\textsubscript{N,k}\textsuperscript{P} & SOCP   & ASM-1~ & PND    & GR-P    \\ 
\hline
1000    & 1.0351                  & 1.0264 & 1.0323 & 1.0398 & 1.0337  \\
2500    & 1.0474                  & 1.0467 & 1.0465 & 1.0438 & -       \\
5000    & -                       & 1.0590 & 1.0540 & 1.0485 & -       \\
10000   & -                       & 1.0659 & 1.0599 & 1.0526 & -       \\
20000   & -                       & 1.0717 & 1.0637 & 1.0573 & -       \\
50000   & -                       & 1.0771 & 1.0702 & -      & -       \\
100000  & -                       & 1.0796 & 1.0727 & -      & -       \\
500000  & -                       & 1.0843 & 1.0776 & -      & -       \\
1000000 & -                       & 1.0855 & 1.0768 & -      & -       \\
\hline
\end{tabular}
\caption{\emph{Objective values for the methods, SOCP, SAA\textsubscript{N,k}\textsuperscript{P}, GR-P, PND and ASM-1. The exact solution to the true problem returns an objective value of 1.0891.}}
\end{subtable}

\end{table}

\begin{table}\ContinuedFloat

\begin{subtable}[h]{1.0\textwidth}

\begin{tabular}{|l|ccccc|} 
\hline
$N$       & SAA\textsubscript{N,k}\textsuperscript{P} & SOCP & ASM-1~ & PND    & GR-P   \\ 
\hline
1000    & 1.2                     & 0.3  & 0.7    & 0.7    & 125.3  \\
2500    & 2632.9                  & 0.5  & 0.8    & 14.2   & $>$3600   \\
5000    & $>$3600                    & 0.7  & 1.1    & 78.0     & -      \\
10000   & -                       & 1.1  & 1.5    & 361.7  & -      \\
20000   & -                       & 1.4  & 2.9    & 1675.3 & -      \\
50000   & -                       & 1.4  & 6.1    & $>$3600   & -      \\
100000  & -                       & 1.3  & 13.4   & -      & -      \\
500000  & -                       & 1.8  & 56.9   & -      & -      \\
1000000 & -                       & 1.8  & 119.2  & -      & -      \\
\hline
\end{tabular}
\caption{\emph{Overall solution times for the methods, SOCP, SAA\textsubscript{N,k}\textsuperscript{P}, GR-P, PND and ASM-1.}}
\end{subtable}

\begin{subtable}[h]{1.0\textwidth}

\begin{tabular}{|l|ccccc|} 
\hline
$N$       & SAA\textsubscript{N,k}\textsuperscript{P} & SOCP & ASM-1~ & PND    & GR-P    \\ 
\hline
1000    & 1.2                     & 0.3  & 0.0181 & 0.0161 & 0.1339  \\
2500    & 2632.9                  & 0.5  & 0.0178 & 0.0135 & -       \\
5000    & -                       & 0.7  & 0.0172 & 0.0136 & -       \\
10000   & -                       & 1.1  & 0.0173 & 0.0138 & -       \\
20000   & -                       & 1.4  & 0.0162 & 0.0140 & -       \\
50000   & -                       & 1.4  & 0.0148 & -      & -       \\
100000  & -                       & 1.3  & 0.0155 & -      & -       \\
500000  & -                       & 1.8  & 0.0119 & -      & -       \\
1000000 & -                       & 1.8  & 0.0120 & -      & -       \\
\hline
\end{tabular}
\caption{\emph{Individual solution times for the methods, SOCP, SAA\textsubscript{N,k}\textsuperscript{P}, GR-P, PND and ASM-1.}}
\end{subtable}

\begin{subtable}[h]{1.0\textwidth}

\begin{tabular}{|l|ccccc|} 
\hline
$N$       & SAA\textsubscript{N,k}\textsuperscript{P} & SOCP & ASM-1~ & PND     & GR-P  \\ 
\hline
1000    & 1                       & 1    & 31.4   & 28.9    & 933   \\
2500    & 1                       & 1    & 34.4   & 461.3   & -     \\
5000    & -                       & 1    & 35.6   & 1588.6  & -     \\
10000   & -                       & 1    & 37.9   & 4390.4  & -     \\
20000   & -                       & 1    & 47.0   & 11421.1 & -     \\
50000   & -                       & 1    & 48.3   & -       & -     \\
100000  & -                       & 1    & 57.1   & -       & -     \\
500000  & -                       & 1    & 62.6   & -       & -     \\
1000000 & -                       & 1    & 61.9   & -       & -     \\
\hline
\end{tabular}
\caption{\emph{Number of times a model is solved for the methods, SOCP, SAA\textsubscript{N,k}\textsuperscript{P}, GR-P, PND and ASM-1.}}
\end{subtable}

\begin{subtable}[h]{1.0\textwidth}

\begin{tabular}{|l|cccccc|} 
\hline
$N$      & $k/N$  & SAA\textsubscript{N,k}\textsuperscript{P} & SOCP  & ASM-1~ & PND   & GR-P   \\ 
\hline
1000    & 0.001 & 0.005                   & 0.001 & 0.005  & 0.018 & 0.005  \\
2500    & 0.010 & 0.011                   & 0.010 & 0.012  & 0.012 & -      \\
5000    & 0.017 & -                       & 0.017 & 0.019  & 0.019 & -      \\
10000   & 0.024 & -                       & 0.024 & 0.025  & 0.026 & -      \\
20000   & 0.030 & -                       & 0.030 & 0.030  & 0.029 & -      \\
50000   & 0.036 & -                       & 0.036 & 0.036  & -     & -      \\
100000  & 0.039 & -                       & 0.039 & 0.039  & -     & -      \\
500000  & 0.045 & -                       & 0.045 & 0.044  & -     & -      \\
1000000 & 0.046 & -                       & 0.047 & 0.046  & -     & -      \\
\hline
\end{tabular}
\caption{\emph{Probability of Violation for the solution of each of the methods, SOCP, SAA\textsubscript{N,k}\textsuperscript{P}, GR-P, PND and ASM-1.}}

\end{subtable}

\begin{subtable}[h]{1.0\textwidth}

\begin{tabular}{|l|cccccc|} 
\hline
$N$       & $k/N$  & SAA\textsubscript{N,k}\textsuperscript{P} & SOCP  & ASM-1~ & PND   & GR-P   \\ 
\hline
1000    & 0.001 & 0.005                   & 0.001 & 0.006  & 0.020 & 0.006  \\
2500    & 0.010 & 0.013                   & 0.011 & 0.014  & 0.014 & -      \\
5000    & 0.017 & -                       & 0.019 & 0.021  & 0.021 & -      \\
10000   & 0.024 & -                       & 0.027 & 0.027  & 0.028 & -      \\
20000   & 0.030 & -                       & 0.032 & 0.032  & 0.031 & -      \\
50000   & 0.036 & -                       & 0.038 & 0.038  & -     & -      \\
100000  & 0.039 & -                       & 0.042 & 0.042  & -     & -      \\
500000  & 0.045 & -                       & 0.048 & 0.047  & -     & -      \\
1000000 & 0.046 & -                       & 0.050 & 0.049  & -     & -      \\
\hline
\end{tabular}
\caption{\emph{Confidence of violation (upper limit) for the solution of each of the methods, SOCP, SAA\textsubscript{N,k}\textsuperscript{P}, GR-P, PND and ASM-1.}}
\end{subtable}
\end{table} 	

\FloatBarrier

\subsection{Discussion of IP Results}

We solve SAA\textsubscript{N,k}\textsuperscript{P} using only 1000 or 2500 scenarios, as the solution time rapidly increases to become impractical, Table 4(b). Using the GR-P method we assume that we cannot define the binding constraints using the slack variables. The example scenario resolves a model $N-1$ times, in each instance with one constraint removed. Each solution takes on average 1.039 seconds. The time itself is not onerous but when one needs $kN$ solutions the method becomes unrealistic.\\

Trialling the PND model we see an order of magnitude time reduction from GR-P for each model solve, Table 4(c). In the trials we test for binding constraints using the slack values of the scenario constraints. Rather than relying on a slack value of zero to define a binding constraint we create the set using a tolerance. In each of the $k$ inner loops of the PND method, we define the set of binding scenario constraints as those with slack values less than a tolerance. In the trials the tolerance is set to $10^{-1}$. The results show a large reduction in each solution time, however the number of solutions required is still of the order $nN$. Thus, the overall time is still impractical.\\

The ASM-1 method gives the best time improvement by only solving a relatively small number of models Table 4(d). The overall time is still practical over all the trials.  It is again clear (Figure 4) that the active set method delivers a better solution in a faster time.\\

\begin{figure}[h] 
\begin{center}
\includegraphics[scale=0.7]{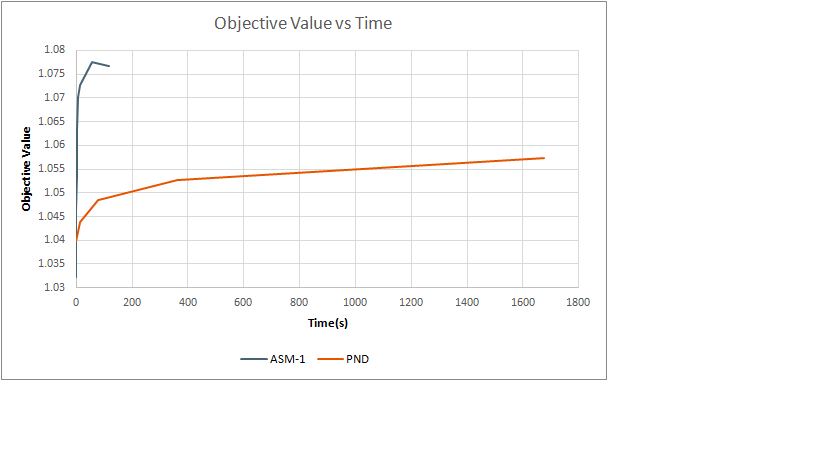}
\end{center}
\vspace{-0.5 in}
\caption{Plot of the Objective Values vs the Time taken to achieve those values for the Active Set and Pool and Discard methods. Results are the mean values of 30 trials. }
\end{figure}

In Table 4(e) we see that the probability of violation closely matches the ratio $k/N$, as expected. In Table 4(f) we see that the $1-\beta$ confidence of the upper limit of the solution is always less than or equal to $\epsilon$. This demonstrates that the solutions obtained are feasible to the original problem AAP as required.\\

Of interest is the behaviour of the solution as $N$ increases. Figures 5(a) to 5(c) show the allocation to each asset in each of the 30 trials for the IP model solved using ASM-1, for the given $N$. The values along the X-axis are the number of trials that allocated zero to that asset. The mean cash allocation is given as CSH.\\

\begin{figure}[h] 
\centering

\begin{subfigure}[h]{1.0\textwidth} 
\begin{center}
\includegraphics[scale=0.6]{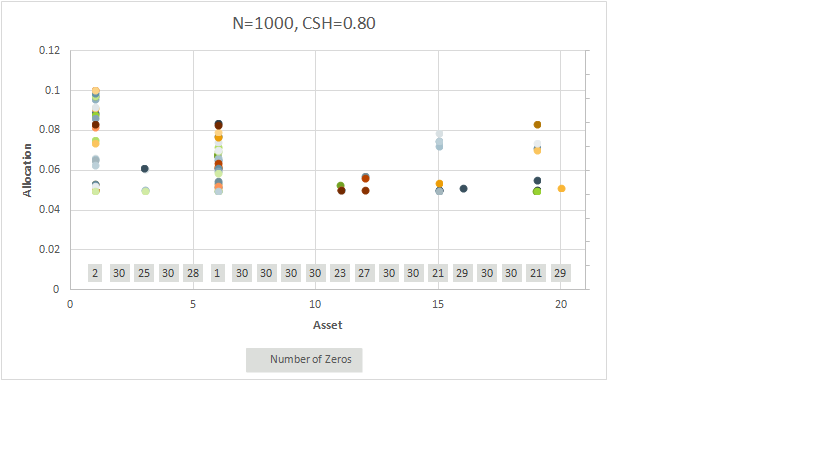}
\end{center}
\vspace{-0.4 in}
\caption{$N=1,000$. Objective values range from 1.0198 to 1.0427; cash allocation ranges from 0.724 to 0.853. }
\end{subfigure}

\begin{subfigure}[h]{1.0\textwidth} 
\begin{center}
\includegraphics[scale=0.6]{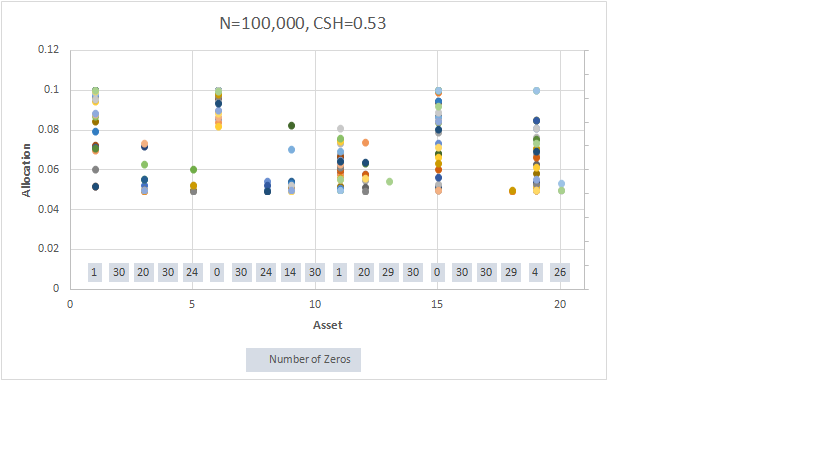}
\end{center}
\vspace{-0.4 in}
\caption{$N=100,000$. Objective values range from 1.0650 to 1.0788; cash allocation ranges from 0.465 to 0.603. }
\end{subfigure}

\begin{subfigure}[h]{1.0\textwidth} 
\begin{center}
\includegraphics[scale=0.6]{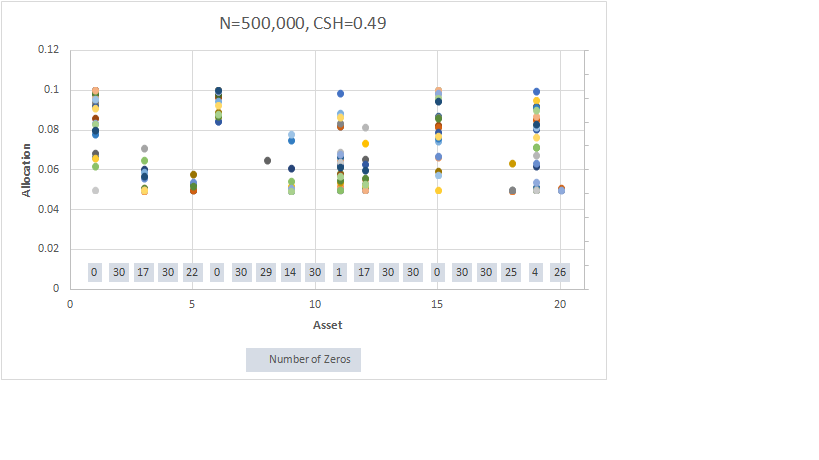}
\end{center}
\vspace{-0.4 in}
\caption{$N=500,000$. Objective values range from 1.0626 to 1.0823; cash allocation ranges from 0.441 to 0.600. }
\end{subfigure}

\caption{Allocation to each of the assets in the 30 trials for $N=1,000$, $N=100,000$ and $N=500,000$.}

\end{figure}

What is clear is that diversification of allocations occurs as the solutions improve. We see that the allocation to cash reduces as the solution improves. We also see that several assets are never considered in any trials. There may be scope to further improve the method by eliminating these assets from consideration. Reducing the number of dimensions $(n)$ will increase the allowable failures $(k)$ in (\ref{eqn3}). Thus there may be some further improvement in the objective values.\\

\section{Conclusion}

What we see is that the Active Set method provides advantages over existing methods used to solve chance-constrained programs with large numbers of constraints. Using the Active Set method, one can find a better solution in a limited amount of time. Reducing the time relative to previous approaches is achieved by solving only a small number of models and keeping their size small. The simple approach of ordering constraints by the value of their dual variables applies to all approaches and is a key contributor to reducing the number of models that need solving.\\

There is no dependence on distribution for the Active Set method. Any shaped distribution can be used to draw samples from, or one can use a-priori data directly to draw samples from. This final point can be useful for some of the sample based optimisation seen in practical applications.\\

The Active Set method is equally applicable to integer-programming models and outperforms all other methods in the examples. This is a significant practical advantage.\\
 
We suggest that further work investigating how to efficiently select the best constraints to add. It is likely that investigating some of the theoretical underpinnings will be fruitful. Expanding the application to other classes of stochastic programming, particularly those considering many scenarios is another avenue to continue this work.\\

\FloatBarrier

\bibliography{sn-bibliography}


\begin{thebibliography}{27}
\ifx \bisbn   \undefined \def \bisbn  #1{ISBN #1}\fi
\ifx \binits  \undefined \def \binits#1{#1}\fi
\ifx \bauthor  \undefined \def \bauthor#1{#1}\fi
\ifx \batitle  \undefined \def \batitle#1{#1}\fi
\ifx \bjtitle  \undefined \def \bjtitle#1{#1}\fi
\ifx \bvolume  \undefined \def \bvolume#1{\textbf{#1}}\fi
\ifx \byear  \undefined \def \byear#1{#1}\fi
\ifx \bissue  \undefined \def \bissue#1{#1}\fi
\ifx \bfpage  \undefined \def \bfpage#1{#1}\fi
\ifx \blpage  \undefined \def \blpage #1{#1}\fi
\ifx \burl  \undefined \def \burl#1{\textsf{#1}}\fi
\ifx \doiurl  \undefined \def \doiurl#1{\url{https://doi.org/#1}}\fi
\ifx \betal  \undefined \def \betal{\textit{et al.}}\fi
\ifx \binstitute  \undefined \def \binstitute#1{#1}\fi
\ifx \binstitutionaled  \undefined \def \binstitutionaled#1{#1}\fi
\ifx \bctitle  \undefined \def \bctitle#1{#1}\fi
\ifx \beditor  \undefined \def \beditor#1{#1}\fi
\ifx \bpublisher  \undefined \def \bpublisher#1{#1}\fi
\ifx \bbtitle  \undefined \def \bbtitle#1{#1}\fi
\ifx \bedition  \undefined \def \bedition#1{#1}\fi
\ifx \bseriesno  \undefined \def \bseriesno#1{#1}\fi
\ifx \blocation  \undefined \def \blocation#1{#1}\fi
\ifx \bsertitle  \undefined \def \bsertitle#1{#1}\fi
\ifx \bsnm \undefined \def \bsnm#1{#1}\fi
\ifx \bsuffix \undefined \def \bsuffix#1{#1}\fi
\ifx \bparticle \undefined \def \bparticle#1{#1}\fi
\ifx \barticle \undefined \def \barticle#1{#1}\fi
\bibcommenthead
\ifx \bconfdate \undefined \def \bconfdate #1{#1}\fi
\ifx \botherref \undefined \def \botherref #1{#1}\fi
\ifx \url \undefined \def \url#1{\textsf{#1}}\fi
\ifx \bchapter \undefined \def \bchapter#1{#1}\fi
\ifx \bbook \undefined \def \bbook#1{#1}\fi
\ifx \bcomment \undefined \def \bcomment#1{#1}\fi
\ifx \oauthor \undefined \def \oauthor#1{#1}\fi
\ifx \citeauthoryear \undefined \def \citeauthoryear#1{#1}\fi
\ifx \endbibitem  \undefined \def \endbibitem {}\fi
\ifx \bconflocation  \undefined \def \bconflocation#1{#1}\fi
\ifx \arxivurl  \undefined \def \arxivurl#1{\textsf{#1}}\fi
\csname PreBibitemsHook\endcsname

\bibitem[\protect\citeauthoryear{Charnes and Cooper}{1959}]{CharnesCooper1959}
\begin{barticle}
\bauthor{\bsnm{Charnes}, \binits{A.}},
\bauthor{\bsnm{Cooper}, \binits{W.W.}}:
\batitle{Chance-constrained programming}.
\bjtitle{Management science}
\bvolume{6}(\bissue{1}),
\bfpage{73}--\blpage{79}
(\byear{1959})
\end{barticle}
\endbibitem

\bibitem[\protect\citeauthoryear{Pagnoncelli et~al.}{2012}]{Pagnoncelli2012}
\begin{barticle}
\bauthor{\bsnm{Pagnoncelli}, \binits{B.K.}},
\bauthor{\bsnm{Reich}, \binits{D.}},
\bauthor{\bsnm{Campi}, \binits{M.C.}}:
\batitle{Risk-return trade-off with the scenario approach in practice: a case
  study in portfolio selection}.
\bjtitle{Journal of Optimization Theory and Applications}
\bvolume{155},
\bfpage{707}--\blpage{722}
(\byear{2012})
\end{barticle}
\endbibitem

\bibitem[\protect\citeauthoryear{{Gurobi Optimization, LLC}}{2023}]{gurobi}
\begin{botherref}
\oauthor{\bsnm{{Gurobi Optimization, LLC}}}:
{Gurobi Optimizer Reference Manual}
(2023).
\url{https://www.gurobi.com}
\end{botherref}
\endbibitem

\bibitem[\protect\citeauthoryear{Jeuken et~al.}{2021}]{Jeuken2021}
\begin{barticle}
\bauthor{\bsnm{Jeuken}, \binits{R.}},
\bauthor{\bsnm{Forbes}, \binits{M.}},
\bauthor{\bsnm{Kearney}, \binits{M.}}:
\batitle{Optimal blending strategies for coking coal using chance constraints}.
\bjtitle{The Journal of the Operational Research Society}
\bvolume{72}(\bissue{12}),
\bfpage{2690}--\blpage{2703}
(\byear{2021})
\end{barticle}
\endbibitem

\bibitem[\protect\citeauthoryear{Meng et~al.}{2020}]{meng2020}
\begin{barticle}
\bauthor{\bsnm{Meng}, \binits{F.}},
\bauthor{\bsnm{Li}, \binits{L.}},
\bauthor{\bsnm{Li}, \binits{T.}},
\bauthor{\bsnm{Fu}, \binits{Q.}}:
\batitle{Optimal allocation model of the water resources in harbin under
  representative concentration pathway scenarios}.
\bjtitle{Water Supply}
\bvolume{20}(\bissue{7}),
\bfpage{2903}--\blpage{2914}
(\byear{2020})
\end{barticle}
\endbibitem

\bibitem[\protect\citeauthoryear{Stuhlmacher and
  Mathieu}{2020}]{stuhlmacher2020}
\begin{barticle}
\bauthor{\bsnm{Stuhlmacher}, \binits{A.}},
\bauthor{\bsnm{Mathieu}, \binits{J.L.}}:
\batitle{Water distribution networks as flexible loads: A chance-constrained
  programming approach}.
\bjtitle{Electric Power Systems Research}
\bvolume{188},
\bfpage{106570}
(\byear{2020})
\end{barticle}
\endbibitem

\bibitem[\protect\citeauthoryear{Wang et~al.}{2020}]{wang2020}
\begin{barticle}
\bauthor{\bsnm{Wang}, \binits{A.}},
\bauthor{\bsnm{Jasour}, \binits{A.}},
\bauthor{\bsnm{Williams}, \binits{B.C.}}:
\batitle{Non-gaussian chance-constrained trajectory planning for autonomous
  vehicles under agent uncertainty}.
\bjtitle{IEEE Robotics and Automation Letters}
\bvolume{5}(\bissue{4}),
\bfpage{6041}--\blpage{6048}
(\byear{2020})
\end{barticle}
\endbibitem

\bibitem[\protect\citeauthoryear{Wu et~al.}{2020}]{wu2020}
\begin{barticle}
\bauthor{\bsnm{Wu}, \binits{P.}},
\bauthor{\bsnm{Chen}, \binits{J.}},
\bauthor{\bsnm{Zhou}, \binits{Z.}},
\bauthor{\bsnm{Zhao}, \binits{Y.}},
\bauthor{\bsnm{Lin}, \binits{F.}}:
\batitle{Determining the optimal location of vehicle inspection facilities
  under uncertainty via new optimization approaches}.
\bjtitle{IEEE Access}
\bvolume{8},
\bfpage{38229}--\blpage{38238}
(\byear{2020})
\end{barticle}
\endbibitem

\bibitem[\protect\citeauthoryear{Zhang et~al.}{2021}]{zhang2021}
\begin{barticle}
\bauthor{\bsnm{Zhang}, \binits{D.}},
\bauthor{\bsnm{Li}, \binits{D.}},
\bauthor{\bsnm{Sun}, \binits{H.}},
\bauthor{\bsnm{Hou}, \binits{L.}}:
\batitle{A vehicle routing problem with distribution uncertainty in deadlines}.
\bjtitle{European Journal of Operational Research}
\bvolume{292}(\bissue{1}),
\bfpage{311}--\blpage{326}
(\byear{2021})
\end{barticle}
\endbibitem

\bibitem[\protect\citeauthoryear{Zhang}{2019}]{Zhang2019}
\begin{barticle}
\bauthor{\bsnm{Zhang}, \binits{P.}}:
\batitle{Multiperiod mean absolute deviation uncertain portfolio selection with
  real constraints}.
\bjtitle{Soft Computing}
\bvolume{23},
\bfpage{5081}--\blpage{5098}
(\byear{2019})
\end{barticle}
\endbibitem

\bibitem[\protect\citeauthoryear{Chowdhury et~al.}{2019}]{Chowdhury2019}
\begin{barticle}
\bauthor{\bsnm{Chowdhury}, \binits{S.}},
\bauthor{\bsnm{Shahvari}, \binits{O.}},
\bauthor{\bsnm{Marufuzzaman}, \binits{M.}},
\bauthor{\bsnm{Francis}, \binits{J.}},
\bauthor{\bsnm{Bian}, \binits{L.}}:
\batitle{Sustainable design of on-demand supply chain network for additive
  manufacturing}.
\bjtitle{Iise Transactions}
\bvolume{51}(\bissue{7}),
\bfpage{744}--\blpage{765}
(\byear{2019})
\end{barticle}
\endbibitem

\bibitem[\protect\citeauthoryear{Wang and Ning}{2018}]{Wang2018}
\begin{barticle}
\bauthor{\bsnm{Wang}, \binits{X.}},
\bauthor{\bsnm{Ning}, \binits{Y.}}:
\batitle{Uncertain chance-constrained programming model for project scheduling
  problem}.
\bjtitle{Journal of the operational research society}
\bvolume{69}(\bissue{3}),
\bfpage{384}--\blpage{391}
(\byear{2018})
\end{barticle}
\endbibitem

\bibitem[\protect\citeauthoryear{Bruninx et~al.}{2017}]{bruninx2017}
\begin{barticle}
\bauthor{\bsnm{Bruninx}, \binits{K.}},
\bauthor{\bsnm{Dvorkin}, \binits{Y.}},
\bauthor{\bsnm{Delarue}, \binits{E.}},
\bauthor{\bsnm{D’haeseleer}, \binits{W.}},
\bauthor{\bsnm{Kirschen}, \binits{D.S.}}:
\batitle{Valuing demand response controllability via chance constrained
  programming}.
\bjtitle{IEEE Transactions on Sustainable Energy}
\bvolume{9}(\bissue{1}),
\bfpage{178}--\blpage{187}
(\byear{2017})
\end{barticle}
\endbibitem

\bibitem[\protect\citeauthoryear{Jiang et~al.}{2017}]{jiang2017}
\begin{barticle}
\bauthor{\bsnm{Jiang}, \binits{Y.}},
\bauthor{\bsnm{Xu}, \binits{J.}},
\bauthor{\bsnm{Shen}, \binits{S.}},
\bauthor{\bsnm{Shi}, \binits{C.}}:
\batitle{Production planning problems with joint service-level guarantee: a
  computational study}.
\bjtitle{International Journal of Production Research}
\bvolume{55}(\bissue{1}),
\bfpage{38}--\blpage{58}
(\byear{2017})
\end{barticle}
\endbibitem

\bibitem[\protect\citeauthoryear{Karimi et~al.}{2021}]{karimi2021}
\begin{barticle}
\bauthor{\bsnm{Karimi}, \binits{H.}},
\bauthor{\bsnm{Ek{\c{s}}io{\u{g}}lu}, \binits{S.D.}},
\bauthor{\bsnm{Carbajales-Dale}, \binits{M.}}:
\batitle{A biobjective chance constrained optimization model to evaluate the
  economic and environmental impacts of biopower supply chains}.
\bjtitle{Annals of Operations Research}
\bvolume{296}(\bissue{1-2}),
\bfpage{95}--\blpage{130}
(\byear{2021})
\end{barticle}
\endbibitem

\bibitem[\protect\citeauthoryear{Pourahmadi et~al.}{2019}]{pourahmadi2019}
\begin{barticle}
\bauthor{\bsnm{Pourahmadi}, \binits{F.}},
\bauthor{\bsnm{Kazempour}, \binits{J.}},
\bauthor{\bsnm{Ordoudis}, \binits{C.}},
\bauthor{\bsnm{Pinson}, \binits{P.}},
\bauthor{\bsnm{Hosseini}, \binits{S.H.}}:
\batitle{Distributionally robust chance-constrained generation expansion
  planning}.
\bjtitle{IEEE Transactions on Power Systems}
\bvolume{35}(\bissue{4}),
\bfpage{2888}--\blpage{2903}
(\byear{2019})
\end{barticle}
\endbibitem

\bibitem[\protect\citeauthoryear{Luedtke and Ahmed}{2008}]{Luedtke2008}
\begin{barticle}
\bauthor{\bsnm{Luedtke}, \binits{J.}},
\bauthor{\bsnm{Ahmed}, \binits{S.}}:
\batitle{A sample approximation approach for optimization with probabilistic
  constraints}.
\bjtitle{SIAM journal on optimization}
\bvolume{19}(\bissue{2}),
\bfpage{674}--\blpage{699}
(\byear{2008})
\end{barticle}
\endbibitem

\bibitem[\protect\citeauthoryear{Campi and Garatti}{2008}]{Campi2008}
\begin{barticle}
\bauthor{\bsnm{Campi}, \binits{M.C.}},
\bauthor{\bsnm{Garatti}, \binits{S.}}:
\batitle{The exact feasibility of randomized solutions of uncertain convex
  programs}.
\bjtitle{SIAM Journal on Optimization}
\bvolume{19}(\bissue{3}),
\bfpage{1211}--\blpage{1230}
(\byear{2008})
\end{barticle}
\endbibitem

\bibitem[\protect\citeauthoryear{Pagnoncelli et~al.}{2009}]{Pagnoncelli2009}
\begin{barticle}
\bauthor{\bsnm{Pagnoncelli}, \binits{B.K.}},
\bauthor{\bsnm{Ahmed}, \binits{S.}},
\bauthor{\bsnm{Shapiro}, \binits{A.}}:
\batitle{Sample average approximation method for chance constrained
  programming: Theory and applications}.
\bjtitle{Journal of optimization theory and applications}
\bvolume{142}(\bissue{2}),
\bfpage{399}--\blpage{416}
(\byear{2009})
\end{barticle}
\endbibitem

\bibitem[\protect\citeauthoryear{Chen et~al.}{2020}]{chen2020}
\begin{barticle}
\bauthor{\bsnm{Chen}, \binits{Y.}},
\bauthor{\bsnm{Li}, \binits{Y.}},
\bauthor{\bsnm{Sun}, \binits{B.}},
\bauthor{\bsnm{Li}, \binits{Y.}},
\bauthor{\bsnm{Zhu}, \binits{H.}},
\bauthor{\bsnm{Chen}, \binits{Z.}}:
\batitle{A chance-constrained programming approach for a zinc hydrometallurgy
  blending problem under uncertainty}.
\bjtitle{Computers \& Chemical Engineering}
\bvolume{140},
\bfpage{106893}
(\byear{2020})
\end{barticle}
\endbibitem

\bibitem[\protect\citeauthoryear{Campi and Garatti}{2011}]{Campi2011}
\begin{barticle}
\bauthor{\bsnm{Campi}, \binits{M.C.}},
\bauthor{\bsnm{Garatti}, \binits{S.}}:
\batitle{A sampling-and-discarding approach to chance-constrained optimization:
  feasibility and optimality}.
\bjtitle{Journal of optimization theory and applications}
\bvolume{148}(\bissue{2}),
\bfpage{257}--\blpage{280}
(\byear{2011})
\end{barticle}
\endbibitem

\bibitem[\protect\citeauthoryear{Kudela and Popela}{2020}]{Kudela2020}
\begin{barticle}
\bauthor{\bsnm{Kudela}, \binits{J.}},
\bauthor{\bsnm{Popela}, \binits{P.}}:
\batitle{Pool \& discard algorithm for chance constrained optimization
  problems}.
\bjtitle{IEEE Access}
\bvolume{8},
\bfpage{79397}--\blpage{79407}
(\byear{2020})
\end{barticle}
\endbibitem

\bibitem[\protect\citeauthoryear{Shapiro et~al.}{2021}]{shapiro2021}
\begin{bbook}
\bauthor{\bsnm{Shapiro}, \binits{A.}},
\bauthor{\bsnm{Dentcheva}, \binits{D.}},
\bauthor{\bsnm{Ruszczynski}, \binits{A.}}:
\bbtitle{Lectures on Stochastic Programming: Modeling and Theory}.
\bpublisher{SIAM}, \blocation{???}
(\byear{2021})
\end{bbook}
\endbibitem

\bibitem[\protect\citeauthoryear{Levin}{1969}]{levin1969}
\begin{barticle}
\bauthor{\bsnm{Levin}, \binits{V.L.}}:
\batitle{Application of e. helly's theorem to convex programming, problems of
  best approximation and related questions}.
\bjtitle{Mathematics of the USSR-Sbornik}
\bvolume{8}(\bissue{2}),
\bfpage{235}
(\byear{1969})
\end{barticle}
\endbibitem

\bibitem[\protect\citeauthoryear{Markowitz}{1952}]{Markowitz1952}
\begin{barticle}
\bauthor{\bsnm{Markowitz}, \binits{H.}}:
\batitle{Portfolio selection}.
\bjtitle{The Journal of Finance}
\bvolume{7}(\bissue{1}),
\bfpage{77}--\blpage{91}
(\byear{1952}).
Accessed 2023-04-17
\end{barticle}
\endbibitem

\bibitem[\protect\citeauthoryear{Xidonas et~al.}{2020}]{Xidonas2020}
\begin{barticle}
\bauthor{\bsnm{Xidonas}, \binits{P.}},
\bauthor{\bsnm{Steuer}, \binits{R.}},
\bauthor{\bsnm{Hassapis}, \binits{C.}}:
\batitle{Robust portfolio optimization: a categorized bibliographic review}.
\bjtitle{Annals of Operations Research}
\bvolume{292}(\bissue{1}),
\bfpage{533}--\blpage{552}
(\byear{2020})
\end{barticle}
\endbibitem

\bibitem[\protect\citeauthoryear{Wallis}{2013}]{wallis2013binomial}
\begin{barticle}
\bauthor{\bsnm{Wallis}, \binits{S.}}:
\batitle{Binomial confidence intervals and contingency tests: mathematical
  fundamentals and the evaluation of alternative methods}.
\bjtitle{Journal of Quantitative Linguistics}
\bvolume{20}(\bissue{3}),
\bfpage{178}--\blpage{208}
(\byear{2013})
\end{barticle}
\endbibitem

\end{thebibliography}


\end{document}